\documentclass[12pt]{article}
\usepackage{setspace}
\usepackage{amssymb}
\usepackage{amsmath}
\usepackage{gdgspace}

\newcommand{\reff}[1]{\ref{#1}}

\numberwithin{equation}{section}

\begin{document}

\vskip .5cm

\begin{center}
\textbf{On the $q$-Analogue of Gamma Functions and Related
Inequalities }
\end{center}
\begin{center}
by
\end{center}
\begin{center}
\textbf{  Taekyun Kim$^{1}$ and  C. Adiga$^{2}$}\\
\vskip .5cm $^{1}$(Home) Ju-Kong APT 103-Dong 1001-Ho Young-Chang
Ri 544-4,\\
Hapcheon-Up Hapcheon-Gun Kyungshang Nam-Do,postal No. 678-802, S.
Korea\\
(Office) Department of Mathematics Education, Kongju National
University, Kongju 314-701,
S. Korea\\
$\text{e-mail:}$$ \text{tkim@kongju.ac.kr / tkim64@hanmail.net}$\\
\vskip .5cm $^{2}$Department of Studies in Mathematics, University
of Mysore, Manasagangotri, Mysore 570006, India\\
$\text{e-mail:c$_{-}$adiga@hotmail.com }$
\end{center}

\vskip .5cm


\begin{center}
\textbf{Abstract}
\end{center}
\begin{spacing}{1.5}


In this paper, we obtain a $q$-analogue of a double inequality
involving the Euler gamma function which was first proved
geometrically by Alsina and Tom\'{a}s [1] and then analytically by
S\'{a}ndor [6].

\begin{tabbing}
\end{tabbing}
\begin{tabbing}
\end{tabbing}


\section{Introduction}
\label{intro} \setcounter{equation}{0} \quad \quad F. H. Jackson
defined the $q$-analogue of the gamma function as
$$\Gamma_q(x)= \frac{(q;q)_\infty}{(q^x;q)_\infty}(1-q)^{1-x},\quad 0<q<1, \text{ cf. [2, 4, 5, 7], }$$

and

$$\Gamma_q(x)= \frac{(q^{-1};q^{-1})_\infty}{(q^{-x};q^{-1})_\infty}
(q-1)^{1-x}q^{x \choose 2},\quad q>1, $$ where
$$(a;q)_\infty = \prod^{\infty}_{n=0}(1-aq^n).$$
It is well known that $\Gamma_q(x)\rightarrow \Gamma(x)$ as
$q\rightarrow 1^{-}$, where $\Gamma(x)$ is the ordinary Euler
gamma function defined by $$\Gamma(x)= \int_{0}^{\infty} e^{-t}
t^{x-1} dt, \quad x>0.$$

Recently Alsina and Tom\'{a}s \cite{1} have proved the following
double inequality on employing a geometrical method:

{\bf Theorem 1.1.} For all $x \in [0,1]$, and for all nonnegative
integers $n$, one has
\begin{equation}
\frac{1}{n!} \leq \frac{\Gamma(1+x)^n}{\Gamma(1+nx)} \leq 1.
\label{eq-1}
\end{equation}
S\'{a}ndor [6] has obtained a generalization of (\reff{eq-1}) by
using certain simple analytical arguments. In fact, he proved that
for all real numbers $a \geq 1$, and all $x \in [0,1]$,

\begin{equation}
\frac{1}{\Gamma(1+a)} \leq \frac{\Gamma(1+x)^a}{\Gamma(1+ax)} \leq
1. \label{eq-2}
\end{equation}
But to prove (\reff{eq-2}), S\'{a}ndor used the following result:

{\bf Theorem 1.2.} For all $x>0$,

\begin{equation}
\frac{\Gamma'(x)}{\Gamma(x)} =-\gamma + (x-1) \sum_{k=0}^{\infty}
\frac{1}{(k+1)(x+k)}. \label{eq-3}
\end{equation}

The main purpose of this paper is to obtain a $q$-analogue
of(\reff{eq-2}). Our proof is simple and straightforward.


\section{MAIN RESULT}

In this section, we prove our main result.

{\bf Theorem 2.1} If $0<q<1, a \geq 1$ and $x \in [0,1]$, then

$$\frac{1}{\Gamma_q(1+a)} \leq
\frac{\Gamma_q(1+x)^a}{\Gamma_q(1+ax)} \leq 1.$$

{\bf Proof.} We have

\begin{equation}
\Gamma_q(1+x) = \frac{(q;q)_{\infty}}{(q^{1+x};q)_{\infty}}
(1-q)^{-x} \label{eq-2.1}
\end{equation}

and

\begin{equation}
\Gamma_q(1+ax) = \frac{(q;q)_{\infty}}{(q^{1+ax};q)_{\infty}}
(1-q)^{-ax}. \label{eq-2.2}
\end{equation}

Taking logarithmic derivatives of (\reff{eq-2.1}) and
(\reff{eq-2.2}), we obtain

\begin{equation}
\frac{d}{dx}\left( \text{log } \Gamma_q(1+x) \right)= -\text{log }
(1-q) + \text{log } q \sum_{n=0}^{\infty}
\frac{q^{1+x+n}}{1-q^{1+x+n}}, \text{ cf. [3, 4, 5], }
\label{eq-2.3}
\end{equation}

and

\begin{equation}
\frac{d}{dx}\left( \text{log }\Gamma_q(1+ax) \right)= -a\text{log
}(1-q) + a\text{log }q \sum_{n=0}^{\infty}
\frac{q^{1+ax+n}}{1-q^{1+ax+n}}.
\label{eq-2.4}
\end{equation}
Since $x\geq 0, a \geq 1, \text{log }q < 0$ and
$$\frac{q^{1+ax+n}}{1-q^{1+ax+n}} - \frac{q^{1+x+n}}{1-q^{1+x+n}}
= \frac{q^{1+ax+n}-q^{1+x+n}}{(1-q^{1+ax+n})(1-q^{1+x+n})} \leq
0,$$ we have
\begin{equation}
\frac{d}{dx}\left( \text{log } \Gamma_q(1+ax) \right) \geq a
\frac{d}{dx}\left( \text{log } \Gamma_q(1+x) \right).
\label{eq-2.5}
\end{equation}

Let
$$g(x)= \text{log } \frac{\Gamma_q(1+x)^{a}}{\Gamma_q(1+ax)},\quad a \geq 1,
x\geq 0.$$
Then
$$g(x)= a \text{log }\Gamma_q(1+x)- \text{log }\Gamma_q(1+ax)$$
and
$$g'(x) = a \frac{d}{dx}\left( \text{log } \Gamma_q(1+x) \right)
- \frac{d}{dx}\left( \text{log } \Gamma_q(1+ax) \right).$$
By (\reff{eq-2.5}), we get $g'(x)\leq 0$, so $g$ is decreasing.
Hence the function
$$f(x) = \frac{\Gamma_q(1+x)^a}{\Gamma_q(1+ax)}, \quad a\geq 1$$
is a decreasing function of $x\geq 0$.
Thus for $x \in [0,1]$ and $a\geq 1$, we have
$$ \frac{\Gamma_q(2)^a}{\Gamma_q(1+a)} \leq
\frac{\Gamma_q(1+x)^a}{\Gamma_q(1+ax)}
  \leq \frac{\Gamma_q(1)^a}{\Gamma_q(1)}. $$
We complete the proof by noting that $\Gamma_q(1)=\Gamma_q(2)=1$.

Remarks. (1) Letting $q$ to 1 in the above theorem. we obtain
(1.2).

(2). Letting $q$ to 1 and then putting $a=n$ in the above theorem,
we get (1.1).


\end{spacing}
\end{document}